\documentclass[12pt]{article}

\usepackage[T2A]{fontenc}
\usepackage[cp1251]{inputenc}
\usepackage[english,russian]{babel}
\usepackage[tbtags]{amsmath}
\usepackage{amsfonts,amssymb}

\usepackage{mathrsfs}

\usepackage[paper]{MZ2Ewin}

\usepackage{graphicx}

\numberwithin{equation}{subsection}

\theoremstyle{plain}

\newcommand{\prodd}{\prod\limits}
\newcommand{\summ}{\sum\limits}
\newcommand{\D}{\mathcal{D}}
\newcommand{\F}{\mathbb{F}}
\renewcommand{\leq}{\leqslant}
\renewcommand{\geq}{\geqslant}

\begin{document}

\udk{517}

\date{00.00.0000}

\author{М.\,Р.~Габдуллин}
\address{Московский государственный университет им. Ломоносова\\
Институт математики и механики УрО РАН.}
\email{Gabdullin.Mikhail@ya.ru}

\title{О квадратах во множестве элементов конечного поля с ограничениями на коэффициенты при разложении по базису}

\markboth{М.\,Р.~Габдуллин}{О квадратах в конечном поле}

\maketitle

\begin{fulltext}
\begin{abstract}

Усилены недавние результаты C.Dartyge, C.Mauduit, A.S\'ark\"ozy в задаче о количестве квадратов среди элементов конечного поля с ограничениями на коэффициенты при разложении по базису. 

\end{abstract}

\footnotetext{Исследование выполнено за счет гранта Российского научного фонда (проект 14-11-00702).}

\subsection{Введение}
\label{subsec1}

При любом фиксированном $b\in\mathbb{N}$, $b\geq2$, каждое число $n\in\mathbb{N}$ единственным образом представимо в системе счисления с основанием $b$:

\begin{equation*}
	n=\sum\limits_{j=0}^{r-1} c_jb^j, \,\,\,\,\,\,\, 0\leq c_j \leq b-1, \,\,\, c_{r-1}\geq 1.
\end{equation*} 
Во многих работах (обширный список приведен в \cite{1} )  изучались арифметические свойства чисел с ''пропущенными'' цифрами, т.е. тех чисел, $b$-ичная запись которых состоит из заданных цифр. 

В \cite{1.5} C.~Dartyge и A.~S\'ark\"ozy рассмотрели аналог этой задачи в конечных полях. Пусть $\mathbb{F}_q$ --- поле из $q=p^r$ элементов, $\{a_1,\ldots,a_r\}$ --- базис $\mathbb{F}_q$ над $\mathbb{F}_p$. Для множества  $\mathcal{D}\subset\mathbb{F}_p$ через $W_{\D}$ будем обозначать множество элементов поля $F_q$, все коэффициенты которых при разложении по базису $\{a_1,\ldots,a_r\}$ принадлежат множеству $\D$. Обозначим через $Q$ множество ненулевых квадратов поля $F_q$. Положим $Q_0=Q\cup\{0\}.$ Будем считать, что $p\geq3$, так как в случае $p=2$ мы имеем $\F_q=Q_0.$

В недавней работе C.~Dartyge, C.~Mauduit, A.~S\'ark\"ozy \cite{1} было показано, что если множество $\D$ достаточно велико, то во множестве $W_{\D}$ имеются квадраты.

\textsc{Теорема A}. \textit{Пусть $\D\subset\F_p$, $2\leq|\D|\leq p-1.$ Тогда}
$$\left||W_{\D}\cap Q_0|-\frac{|W_{\D}|}{2}\right| \leq \frac{1}{2\sqrt{q}}\left(|\D|+p\sqrt{p-|\D|}\right)^r.
$$
Эта оценка нетривиальна, если $|\D|\geq\frac{(\sqrt5-1)p}{2}(1+o_p(1))$.

В случае, когда множество $\D$ состоит из последовательных чисел, в этой же работе был получен аналог предыдущей теоремы.

\textsc{Теорема B}. \textit{Пусть $\D=\{0,\ldots,t-1\}$, где $2\leq t\leq p-1$. Тогда}
$$\left||W_{\D}\cap Q_0|-\frac{|W_{\D}|}{2}\right| \leq \frac12\left(C(p,t)t\sqrt{p}\right)^r,
$$
\textit{где} 
$$ C(p,t)=\begin{cases} 
\frac{\log p}{t}+\frac1t\left(\frac43-\frac{\log 3}{2}\right)+\frac1p,&\text{если $2\leq t<p-2$,} \\
\frac2p+\frac{2}{\pi(p-1)}(1-\log(2\sin\frac{\pi}{2p})),&\text{если $t=p-2$}. \end{cases} 
$$
Эта оценка нетривиальна, если $t\gg \sqrt{p}\log p$.

В настоящей работе будут доказаны следующие две оценки на количество квадратов во множестве $W_{\D}$, из которых вытекает существование квадратов при ограничениях на размер множества $\D$ более слабых, чем в теореме А.

\textbf{Теорема 1.} \textit{Пусть} $2r-1\leq p^{1/2}$. \textit{Тогда справедлива оценка}
\begin{equation*}
\left| |W_{\D}\cap Q| - \frac{|W_{\D}|}{2}\right| \leq
\frac12|\D|^{1/2}\left(p^{1/4}(2r-1)^{1/2}|\D|^{r-1} + \frac14p^{3/4}r^{3/2}+p^{1/2}\right)+\frac12 \, .
\end{equation*}	
\textit{В частности, если} $\delta=\left(\sqrt p(2r-1)\right)^{2-r}$ \textit{и} $|\D|\geq(1+\delta)(2r-1)p^{1/2}$, \textit{то} $|W_{\D}\cap Q|\geq1.$

\bigskip

\textbf{Теорема 2.} \textit{При любых натуральных} $\nu$ и $1\leq k \leq r-1$ \textit{справедлива оценка}
\begin{equation*}
\left| |W_{\D}\cap Q| - \frac{|W_{\D}|}{2}\right|< \frac12|\D|^{(r-k)(1-1/2\nu)}\left((2\nu)^{\nu}|\D|^{k\nu}q+|\D|^{2k\nu}4\nu q^{1/2}\right)^{1/2\nu}+\frac12.
\end{equation*}	
\textit{Кроме того, если} $r\geq20,$  $C(r)=\exp\left(\frac{4\log r+8}{r}\right)=1+o(1), \, r\to\infty$, \textit{то при} $|\D|\geq  C(r) p^{\frac12}\exp\left(\frac{\log p+4\log\log p}{r}\right)$ \textit{имеем} $|W_{\D}\cap Q|\geq1$.

\bigskip

В частности, из теоремы 2 следует, что при большом $r$ во множестве $W_{\D}$ есть квадраты уже при $|\D|>p^{1/2}$. Отметим, что при $r \gg \frac{\log p}{\log\log p}$, более точный результат дает теорема 2, а иначе --- теорема 1.

При малых r теорему B также можно усилить, пользуясь оценкой сумм характеров, полученной в работе С.\,В.~Конягина \cite{2}.

\textsc{Теорема C}. \textit{Пусть $\varepsilon\in(0,1/4]$, $\chi$ -- нетривиальный мультипликативный характер в $\F_q$, $N_i$, $H_i$ -- целые числа, $p^{1/4+\varepsilon}\leq H_i \leq p$, $i=1,\ldots,r$,  и }
\begin{equation}
B=\left\{\sum\limits_{i=1}^{r} x_ia_i : N_i+1 \leq x_i \leq N_i+H_i, \quad i=1,\ldots,r\right\}. \label{B}
\end{equation}
\textit{Тогда}
$$\left|\sum\limits_{x\in B}\chi(x)\right| \ll \frac{r^{O(1)}}{\varepsilon} p^{-\varepsilon^2/2}|B|.
$$

Рассуждая стандартным образом (см., например, начало доказательства теоремы~1), из теоремы С нетрудно вывести следующий результат.

\textsc{Следствие.} \textit{Пусть $\D=\{0,1,\ldots,t-1\}$, $\varepsilon>0$, $t\geq p^{1/4+\varepsilon}$. Тогда справедлива оценка}
$$\left| |W_{\D}\cap Q| - \frac{|W_{\D}|}{2}\right| \ll   \frac{r^{O(1)}}{\varepsilon} p^{-\varepsilon^2/2}|W_{\D}|.
$$
\textit{В частности, если} $\varepsilon \geq C\left(\sqrt{\frac{\log r}{\log p}}+\frac{\log\log p}{(\log p)^{1/2}(\log\log p + \log r)^{1/2}}\right)$ \textit{с некоторой абсолютной постоянной $C>0$, то} $|W_{\D}\cap Q|\geq1$.

\bigskip

После того, как данная работа была подана в печать, в открытом доступе появилась работа R.~Dietmann, C.~Elsholtz, I.\,E.~Shparlinski \cite{2.5}, в которой была рассмотрена более общая задача. Пусть $D_1,\ldots,D_r$ -- подмножества  $\mathbb{F}_p$. Положим
$$ W=W(D_1,\ldots,D_r)=\left\{ x_1a_1+\ldots+x_ra_r \,|\, x_i\in D_i  \right\}. 
$$
Авторы работы \cite{2.5} отмечают, что доказательство теоремы А \cite{1} переносится на случай, когда множества $D_i$ различны, а именно, при $\min\limits_{1\leq i\leq r} |D_i| \geq \frac{(\sqrt5-1)p}{2}(1+o_p(1))$ справедливо $|W\cap Q_0|\geq 1$, и доказывают более сильное утверждение.  

\textsc{Теорема} (\cite{2.5}, теорема 3.5). \textit{Для любого $\varepsilon>0$ существует $\delta>0$ такое, что для любых множеств $D_1,\ldots,D_r$, удовлетворяющих условиям}
$$ \prod_{i=1}^r |D_i| \geq p^{(1/2+\varepsilon)r^2/(r-1)}
$$
\textit{и}
$$ \min\limits_{1\leq i\leq r} |D_i| \geq p^{\varepsilon}
$$
\textit{справедливо} $|W\cap Q_0|=\left(\frac12+O(p^{-\delta})\right)|W|$.

По аналогии с работой \cite{2.5}, теорема B также может быть перенесена на случай различных множеств $D_i$ (см. \cite{2.0}).

В разделе 2 мы приводим необходимые вспомогательные результаты. В разделах 3 и 4 приводятся доказательства теорем 1 и 2 соответственно. Далее, в работе \cite{2} множитель $\frac{r^{O(1)}}{\varepsilon}$ не был выписан явно; для полноты мы докажем теорему С в сформулированном виде в разделе 5.

\bigskip

Автор благодарен С.\,В.~Конягину за постановку задачи и внимание к работе.

\subsection{Вспомогательные результаты}
\label{subsec2}

Приведём леммы, которые понадобятся \newline нам в дальнейшем.

\bigskip

\textsc{Лемма D \cite{3},\cite{3.5}}. \textit{Пусть $\chi$ -- мультипликативный характер порядка $s$ в $\F_{q}$ и $\alpha,\beta\in\F_{q}$ -- несопряжённые порождающие элементы $\F_{q}$ над $\F_p$. Тогда}

$$\left|\sum\limits_{\xi\in\F_p}\chi\left((\xi+\alpha)(\xi+\beta)^{s-1}\right)\right| \leq (2r-1)p^{1/2}.
$$

\bigskip

\textsc{Лемма E \cite{4}}. \textit{Пусть $t$ -- целое число, $1\leq t<q$, \, $\chi_1,\ldots,\chi_t$ -- мультипликативные характеры в $F_q$, причём для некоторого $i$ $\chi_i \neq \chi_0$, где $\chi_0$ -- главный характер. Пусть, далее, $h_1,\ldots,h_t$ -- различные элементы $F_q$ и}

$$ S=\sum\limits_{a\in F_q} \chi_1(a+h_1)\chi_2(a+h_2)\cdots\chi_t(a+h_t).
$$
\textit{Тогда}
$$ |S|\leq (t-t_0-1)q^{1/2}+t_0+1
$$
\textit{где $t_0$ -- число характеров $\chi_i$, для которых $\chi_i=\chi_0.$}

\bigskip

Хорошо известны оценки сумм характеров по суммам множеств (см., например, \cite{5}, лемма 2). Нам будет удобно использовать оценку следующего вида.

\textsc{Лемма 1.} \textit{Для любых $\nu\in\mathbb{N}$, $U,V\subset\F_q$ и квадратичного характера $\chi$ на $\F_q$ справедливо}

$$ \left|\sum\limits_{u\in U}\sum\limits_{v\in V} \chi(u+v)\right| \leq |U|^{1-1/2\nu}\left(\frac{(2\nu)!}{\nu!}|V|^{\nu}q+|V|^{2\nu}4\nu q^{1/2}\right)^{1/2\nu}.
$$

\textit{Док--во.} Имеем

$$ \left|\sum\limits_{u\in U}\sum\limits_{v\in V} \chi(u+v)\right| \leq \sum\limits_{u\in U}\left|\sum\limits_{v\in V} \chi(u+v)\right|\leq
|U|^{1-1/2\nu}S^{1/2\nu},
$$
где
$$ S=\sum\limits_{u\in U}\left| \sum\limits_{v\in V} \chi(u+v)\right|^{2\nu} \leq \sum\limits_{(v_1,\ldots,v_{2\nu})\in V^{2\nu}}\left|\sum\limits_{a\in\F_q}\chi(a+v_1)\cdots\chi(a+v_{2\nu})\right|.
$$
Последнюю сумму разобьем на две: на сумму $S_1$ по тем наборам $(v_1,\ldots,v_{2\nu}),$ в которых значение каждой компоненты встречается чётное число раз, и на сумму $S_2$ по всем остальным наборам. Ясно, что $S_1 \leq C_{2\nu}^{\nu}|V|^{\nu}\nu!q=  \frac{(2\nu)!}{\nu!}|V|^{\nu}q$. Для оценки внутренней суммы для наборов из суммы $S_2$ воспользуемся леммой E; каждая полученная оценка будет не больше, чем $(2\nu-1)q^{1/2}+2\nu  <  4\nu q^{1/2}$. Поэтому
$$
S\leq S_1+S_2 < \frac{(2\nu)!}{\nu!}|V|^{\nu}q+|V|^{2\nu}4\nu q^{1/2}.
$$
Лемма доказана.

\subsection{Доказательство теоремы 1}
\label{subsec3}

Через $\chi$ обозначим квадратичный характер на $\mathbb{F}_q$; cчитаем, что $\chi(0)=0$. Пусть $0$ не принадлежит $\D$. Тогда

$$|W_{\D}\cap Q|=\frac12\sum\limits_{x\in W_{\D}} (1+\chi(x))=\frac12|W_{\D}|+\frac12\sum\limits_{x\in W_{\D}} \chi(x)=\frac12|\D|^r+\frac12\sum\limits_{x\in W_{\D}} \chi(x).
$$
Если же $0\in \D,$ то
$$|W_{\D}\cap Q|=\frac12\sum\limits_{x\in W_{\D}\setminus\{0\}} (1+\chi(x))=\frac12(|\D|^r-1)+\frac12\sum\limits_{x\in W_{\D}} \chi(x).
$$
Таким образом, всегда справедлива оценка
\begin{equation}
\left| |W_{\D}\cap Q| - \frac{|W_{\D}|}{2}\right| \leq \frac12\left|\summ_{x\in W{\D}}\right|+\frac12, \label{est1}
\end{equation}
и нужно для доказательства нужно оценить сумму характеров. Положим $b_j=a_j/a_1$. Тогда $b_1=1$ и $\{1,b_2,\ldots,b_r\}$ --- базис. Имеем
\begin{equation}
\left| \sum\limits_{x\in W_{\D}} \chi(x)\right| \leq \sum\limits_{c_1\in \D}\left| \sum\limits_{c_2,\ldots,c_r\in \D} \chi(c_1a_1+\ldots+c_ra_r)\right| \leq |\D|^{1/2}A^{1/2}, \label{1}
\end{equation}
где
\begin{multline*} A=\sum\limits_{c_1\in \D}\left| \sum\limits_{c_2,\ldots,c_r\in \D} \chi(c_1a_1+\ldots+c_ra_r)\right|^2 =\\
\sum\limits_{c_1\in \D}\left| \sum\limits_{(c_2,\ldots,c_r)\in \D^{r-1}} \chi(c_1+c_2b_2\ldots+c_rb_r)\right|^2.
\end{multline*}
Пусть $\D_d$ --- множество тех наборов $(c_2,\ldots,c_r)\in \D^{r-1}$, для которых элемент $c_2b_2+\ldots+c_rb_r$ лежит в подполе порядка $p^d$ и не лежит ни в каком подполе меньшего порядка. Ясно, что $\D^{r-1}=\bigsqcup\limits_{d|r} \D_d$, причем $\D_1=\{0\}$, если $0\in \D$, и $\D_1=\emptyset$ иначе. Для $d|r$ определим функцию $f_d(x) \colon \D\to\mathbb{C}$,  $f_d(c)=\sum\limits_{(c_2,\ldots,c_r)\in \D_d} \chi(c+c_2b_2+\ldots+c_rb_r)$. Напомним, что $l_2$-норма функции $g \colon \D\to\mathbb{C}$ определяется как $\| g \|_2 = \left(\sum\limits_{x\in\D}|g(x)|^2 \right)^{1/2}$. Тогда в силу неравенства треугольника

\begin{equation}A^{1/2}=\left\| \sum_{d|r} f_d \right\|_2 \leq \sum\limits_{d|r} \|f_d\|_2 = \sum\limits_{d|r} A_d^{1/2}. \label{A^{1/2}}
\end{equation}
где
$$ A_d= \sum\limits_{x\in \D}\left| \sum\limits_{(c_2,\ldots,c_r)\in \D_d} \chi(x+c_2b_2+\ldots+c_rb_r)\right|^2 .
$$
По определению множества $\D_d$ при любом $(c_2,\ldots,c_r)\in \D_d$ элемент $c_2b_2+\ldots+c_rb_r$ порождает подполе порядка $p^d$. Учитывая, что каждый такой элемент имеет не более $d$ сопряженных, и применяя лемму D к парам несопряженных элементов, при $d>1$ имеем
\begin{multline*} A_d\leq \sum\limits_{(c_2,\ldots,c_r), (c'_2,\ldots,c'_r)\in \D_d} \left|\sum\limits_{x\in F_p}  \chi(x+c_2b_2+\ldots+c_rb_r)\overline{\chi}(x+c'_2b_2+\ldots+c'_rb_r)\right|
\leq \\
\sum\limits_{(c_2,\ldots,c_r)\in \D_d} \left(dp+(|\D_d|-d)(2d-1)p^{1/2}\right)\leq(2d-1)p^{1/2}|\D_d|^2+dp|\D_d|.
\end{multline*}
Кроме того, $A_1\leq |\D|\leq p.$ Обозначим $\mathcal{J}=\{d|r : d>1 \,\mbox{и}\, \D_d \neq\emptyset \}$. Тогда при $d\in\mathcal{J}$ в силу неравенства $\sqrt{A+B}\leq \sqrt{A}\left(1+\frac{B}{2A}\right)$, верного при всех положительных $A,B$, получаем
$$ A_d^{1/2}\leq (2d-1)^{1/2}p^{1/4}|\D_d|+\frac{dp^{3/4}}{2(2d-1)^{1/2}}.
$$
Из этой оценки и неравенства (\ref{A^{1/2}}) имеем
$$ A^{1/2}\leq p^{1/4}S_1+\frac12p^{3/4}S_2+p^{1/2},
$$
где
\begin{gather*}
S_1=\sum\limits_{d\in \mathcal{J}} (2d-1)^{1/2}|\D_d|,\quad
S_2=\sum\limits_{d\in\mathcal{J}}\frac{d}{(2d-1)^{1/2}}.
\end{gather*}
Учитывая, что $\sum\limits_{d|r} |\D_d|=|\D|^{r-1}$, получаем
$$
S_1\leq(2r-1)^{1/2}|\D|^{r-1},\quad
S_2\leq\sum\limits_{d\in \mathcal{J}} d^{1/2} \leq \frac12r^{3/2}.
$$
(Последняя оценка проверяется непосредственно при $2\leq r\leq7$, а при $r \geq 8$ вытекает из неравенств $r^{1/2}\leq\frac16r^{3/2}$ и $\sum\limits_{d\leq r/2} d^{1/2}\leq \frac23(r/2+1)^{3/2}\leq \frac13r^{3/2}.$)
Значит,
$$A^{1/2}\leq p^{1/4}(2r-1)^{1/2}|\D|^{r-1}+\frac14p^{3/4}r^{3/2}+p^{1/2}.
$$
Подставляя последнее неравенство в (\ref{1}), получим
$$
\left| \sum\limits_{x\in W_{\D}} \chi(x)\right| \leq |\D|^{1/2}\left(p^{1/4}(2r-1)^{1/2}|\D|^{r-1} + \frac14p^{3/4}r^{3/2}+p^{1/2}\right).
$$
Отсюда и из (\ref{est1}) вытекает первое утверждение теоремы. Далее, во множестве $W_{\D}$ есть квадраты, если правая часть последнего неравенства $<|\D|^r-1$. Это равносильно условию
$$|\D|^{r-1}\left(|\D|^{1/2}-(2r-1)^{1/2}p^{1/4}\right) > \frac14p^{3/4}r^{3/2}+p^{1/2}+|\D|^{-1/2}.
$$
Покажем теперь, что последнее неравенство выполнено при $|\D|\geq (1+\delta)(2r-1)p^{1/2}$, где $\delta=\left(\sqrt p(2r-1)\right)^{2-r}$. В силу того, что $\sqrt{1+\delta}-1\geq\frac{\delta}{2\sqrt2}$ при $\delta\in(0,1]$, а также $2r-1\geq\frac32r$ при $r\geq2$, имеем
\begin{multline*}
|\D|^{r-1}\left(|\D|^{1/2}-(2r-1)^{1/2}p^{1/4}\right) \geq (1+\delta)^{r-1}(2r-1)^{r-1/2}(p^{1/2})^{r-1/2}\frac{\delta}{2\sqrt2}= \\
\frac{(1+\delta)^{r-1}}{2\sqrt2}(2r-1)^{3/2}p^{3/4}\geq\frac{3\sqrt3}{8}p^{3/4}r^{3/2}>\frac14p^{3/4}r^{3/2}+p^{1/2}+|\D|^{-1/2}.
\end{multline*}
Теорема доказана.

\subsection{Доказательство теоремы 2}
\label{subsec4}

Введём натуральные параметры $k$ и $\nu$, которые выберем позже. Положим
$$ U=\left\{\sum\limits_{j=1}^{r-k} c_ja_j : c_j\in \D\right\}, \quad V=\left\{\sum\limits_{j=r-k+1}^{r}c_ja_j : c_j\in \D\right\} ,
$$
где $1\leq k\leq r-1$. Тогда $W_{\D}=U+V$ и по лемме 1
\begin{equation}
\left|\sum\limits_{x\in W_D}\chi(x)\right| < |\D|^{(r-k)(1-1/2\nu)}\left((2\nu)^{\nu}|\D|^{k\nu}q+|\D|^{2k\nu}4\nu q^{1/2}\right)^{1/2\nu}.  \label{2}
\end{equation}
Отсюда и из (\ref{est1}) вытекает первое утверждение теоремы. Далее, положим $a=|\D|^{-r}$. Во множестве $W_{\D}$ есть квадраты, если оценка (\ref{2}) нетривиальна, т.е. если правая часть $\leq|\D|^r-1$. Это равносильно условию
$$(2\nu)^{\nu}|\D|^{k\nu}q+|\D|^{2k\nu}4\nu q^{1/2} \leq |\D|^{2k\nu+r-k}(1-a)^{2\nu},
$$
или
$$|\D|^{k\nu}\left(|\D|^{r-k}(1-a)^{2\nu}-4\nu q^{1/2}\right) \geq (2\nu)^{\nu}q.
$$
Пусть $|\D|^{r-k}\geq 5\nu q^{1/2}(1-a)^{-2\nu}$, или
\begin{multline}
|\D|\geq (5\nu)^{\frac{1}{r-k}}p^{\frac{r}{2(r-k)}}(1-a)^{-\frac{2\nu}{r-k}}
=(5\nu)^{\frac{1}{r-k}}p^{\frac12+\frac{k}{2(r-k)}}(1-a)^{-\frac{2\nu}{r-k}}=\\ p^{1/2}\exp\left(\frac{1}{r-k}(\log 5\nu +\frac12k\log p + 2\nu\log(1-a)^{-1})\right)
\label{D}.
\end{multline}
Тогда $|\D|^k \geq \left(5\nu q^{1/2}\right)^{\frac{k}{r-k}}(1-a)^{-\frac{2\nu k}{r-k}}$ и
$$|\D|^{k\nu}\left(|\D|^{r-k}(1-a)^{2\nu}-4\nu q^{1/2}\right) \geq  \nu q^{1/2}\left(5\nu q^{1/2}\right)^{\frac{k\nu}{r-k}}(1-a)^{-\frac{2k\nu^2}{r-k}}.
$$
Значит, во множестве $W_{\D}$ имеются квадраты, если выполнено условие (\ref{D}) и неравенство
\begin{equation}
\nu \left(5\nu q^{1/2}\right)^{\frac{k\nu}{r-k}}(1-a)^{-\frac{2k\nu^2}{r-k}} \geq (2\nu)^{\nu}q^{1/2}.
\label{3}
\end{equation}
Вместо (\ref{3}) потребуем условие
\begin{equation}\left(q^{1/2}\right)^{\frac{k\nu}{r-k}-1} \geq  (2\nu)^{\nu} \label{cond}
\end{equation}
(Так как $a\in(0,1)$, то это условие является более сильным, чем (\ref{3})). Запишем (\ref{cond}) в виде
$$\exp\left(\frac12r\log p (k\nu+k-r)\right) \geq \exp\left((r-k)\nu\log 2\nu \right),
$$
то есть
$$r\log p (k\nu+k-r) \geq 2\nu(r-k)\log 2\nu,
$$
или
$$\log p \left(k+\frac{k}{\nu}-\frac{r}{\nu}\right) \geq 2\left(1-\frac{k}{r}\right)\log 2\nu.
$$
Последнее выполнено, если
$$\log p \left(k-\frac{r}{\nu}\right) \geq 2\log 2\nu ,
$$
т.е. если
\begin{equation} k\log p \geq 2\log 2\nu+\frac{r}{\nu}\log p. \label{k}
\end{equation}

Таким образом, если выполнены (\ref{D}) и (\ref{k}), то $|W_{\D}\cap Q|\geq1$.

В наших интересах выбрать $k$ как можно меньше. Минимум в правой части (\ref{k}) достигается при $\nu=\frac12 r\log p$. Положим $\nu=\left[\frac12 r\log p\right]$. Тогда для выполнения (\ref{k}) достаточно брать $k\geq (\log p)^{-1}(2\log r+2\log\log p+4).$
Положим
$$ k=\left[ \frac{2\log r+2\log\log p+4}{\log p} \right]+1.
$$
Так как при $u\geq3$ функция $\frac{2\log\log u+4}{\log u}$ убывает, то $k\leq \left[ \frac{2\log r+2\log\log 3+4}{\log 3} \right]+1 \leq r/2$ при $r\geq20$. Далее, имеем $\log (1-a)^{-1} \leq \frac{a}{1-a}\leq \frac43a = \frac43|\D|^{-r}\leq \frac43p^{-r/2}$. Поэтому при выбранных $\nu$ и $k$ показатель экспоненты в последней строке (\ref{D}) будет не больше, чем
\begin{multline*}
\frac2r\left(\log(5/2)+ \log (r\log p) + \log r+\log\log p + 2 +\frac12\log p+\frac43p^{-r/2}r\log p  \right) = \\
\frac2r\left(2\log r+\frac12\log p+ 2\log\log p+2+\log(5/2)+\frac43p^{-r/2}r\log p \right) \leq \\
 \frac1r\left(4\log r+\log p+ 4\log\log p+8 \right) .
\end{multline*}
Отсюда вытекает второе утверждение теоремы.

\subsection{Доказательство теоремы C}
\label{subsec5}

В этом разделе мы полностью следуем работе [2], лишь отслеживая зависимость констант от параметров $r$ и $\varepsilon$.

\smallskip

\textsc{Лемма 2}. \textit{Пусть множество $B$ определено соотношением (\ref{B}), где $N_i,H_i>0$ -- целые числа, причем $H_1=\ldots=H_r \leq p^{1/2}$. Тогда уравнение}
\begin{equation}
x^1x^2=x^3x^4, \,\,\, x^1,x^2,x^3,x^4 \in B,   \label{eq}
\end{equation}
\textit{имеет не более $r^{O(r)}|B|^2\log p$ решений}.

\smallskip

\textit{Док-во}. Положим $H=H_1=\ldots=H_r,$

$$ Z=\frac{B\setminus\{0\}}{B\setminus\{0\}}=\{z\in \F_q : \exists x,y\in B\setminus\{0\} \, xz=y\}.
$$
Если $x^1,x^2,x^3,x^4 \in B$, $x^1x^2=x^3x^4$, и $(x^1,x^4)\neq(0,0)$, $(x^2,x^3)\neq(0,0)$, то для некоторого $z\in Z$ верно $x^1z=x^3$, $x^4z=x^2$. Поэтому число $E$ решений уравнения (\ref{eq}) оценивается следующим образом:

$$ E \leq 2|B|^2 + \sum\limits_{z\in Z} f^2(z),
$$
где $f(z)$ -- число решений уравнения
$$xz=y, \quad\quad x,y\in B.
$$
Обозначим
$$B_0=\{x_1a_1+\ldots+x_ra_r : -H\leq x_j\leq H\},
$$
$$Z_0=\frac{B_0\setminus\{0\}}{B_0\setminus\{0\}}, \quad\quad f_0(z)=|\{ (x,y)\in B^2_0 : xz=y \}|.
$$
Заметим, что $f(z)\leq f_0(z)$ и $f_0(z)=1$ при $z\in \F_q^{*}\setminus Z_0$. Поэтому
$$\sum\limits_{z\in Z} f^2(z) \leq \sum\limits_{z\in Z_0} f_0^2(z) +\left|Z\setminus Z_0 \right|.
$$
Так как $|Z|\leq |B|^2$, получаем
\begin{equation} E\leq 3|B|^2 + \sum\limits_{z\in Z_0} f_0^2(z), \label{4}
\end{equation}
и остается оценить последнюю сумму.

Для фиксированного $z\in Z_0$ определим решетку $\Gamma$ в $\mathbb{Z}^{2r}$:

$$ \Gamma = \Gamma_z = \left\{ (x_1,\ldots,x_r,y_1,\ldots,y_r)\in \mathbb{Z}^{2r} : z\sum\limits_{i=1}^r x_ia_i = \sum\limits_{i=1}^r y_ia_i \right\}.
$$
Для фиксированных $x_1,\ldots,x_r \in \mathbb{Z}$ условие $(x_1,\ldots,x_r,y_1,\ldots,y_r)\in \Gamma$ определяет вычет по модулю $p$ каждого из чисел $y_1,\ldots,y_r$. Значит, количество точек решетки $\Gamma$ в большом кубе имеет следующую асимптотику при $M\to\infty$:

$$ \left| \{ (x_1,\ldots,y_r)\in \Gamma : |x_i|\leq M, |y_i|\leq M, i=1,\ldots,r\}\right| = \frac{(2M)^{2r}}{p^r}(1+o(1)).
$$
Поэтому
\begin{equation} \mbox{mes}(\mathbb{R}^{2r}/\Gamma) = p^r.  \label{5}
\end{equation}
Определим куб $D \subset \mathbb{R}^{2r}$:

$$ D = \{ (x_1,\ldots,y_r)\in\mathbb{R}^{2r} : |x_i|\leq H, |y_i|\leq H, i=1,\ldots,r\}.
$$
Отметим, что
\begin{equation} f_0(z)=\left| D\cap\Gamma_z \right| . \label{6}
\end{equation}
Напомним, что \textit{$i$-й последовательный минимум}
$$\lambda_i = \lambda_i(z) = \lambda_i(D,\Gamma_z)
$$
множества $D$ по отношению к $\Gamma_z$ определяется как минимальное число $\lambda$, при котором множество $\lambda D$ содержит $i$ линейно независимых векторов решетки $\Gamma_z$, $i=1,\ldots,2r$. Очевидно, $\lambda_1(z)\leq\ldots\leq\lambda_{2r}(z)$, причем условие $\lambda_1(z)\leq 1$ равносильно тому, что $z\in Z_0$. Вторая теорема Минковского (см., например, [6], теорема 3.30) утверждает, что

$$ \frac{2^{2r}}{(2r)!} \leq \frac{\lambda_1\ldots\lambda_{2r}\mbox{mes}D}{\mbox{mes}(\mathbb{R}^{2r}/\Gamma)} .
$$
Учитывая (\ref{5}), получим
\begin{equation} \lambda_1\ldots\lambda_{2r} \geq \frac{p^r}{(2r)!H^{2r}} . \label{7}
\end{equation}
Хорошо известно (см. [7], предложение 2.1), что число $f_0(z)$ точек решетки $\Gamma$ в кубе $D$ удовлетворяет неравенству
\begin{multline} f_0(z) \leq \prodd_{i=1}^{2r} \left(\frac{2i}{\lambda_i}+1\right) \leq  2^{2r}(2r)!\prodd_{i=1}^{2r}\left(\frac{1}{\lambda_i}+\frac{1}{2i}\right)\leq
r^{O(r)}\prodd_{i=1}^{2r} \max\left(\frac{1}{\lambda_i},1\right). \label{8}
\end{multline}

Определим полярную решетку $\Gamma^{*}=\Gamma^{*}_z$ как множество векторов $(u_1,\ldots,u_{2r})\in \mathbb{R}^{2r}$ таких, что
$$ \summ_{i=1}^{r} u_ix_i +\summ_{i=1}^r u_{r+i}y_i \in \mathbb{Z}
$$
для всех $(x_1,\ldots,x_r,y_1,\ldots,y_r)\in\Gamma$. Заметим, что $\Gamma^*_z \subset p^{-1}\mathbb{Z}^{2r}$, так как $p\mathbb{Z}^{2r}\subset \Gamma_z$. Обозначим через $\lambda_1^*=\lambda_1^*(z)$ первый последовательный минимум решетки $\Gamma_z^*$ по отношению к множеству
$$ D^*=\left\{(u_1,\ldots,u_{2r}) : \summ_{i=1}^{2r} |u_i| \leq \frac{1}{H} \right\}
$$
Согласно [8] (см. предложение 3.6), имеем
\begin{equation} \lambda_{2r}\lambda_1^* \ll r(1+\log r)^{1/2} = r^{O(1)}. \label{9}
\end{equation}
Обозначим
$$
\nu=\nu(z)=\min\left(\lambda_1H,\frac{p\lambda_1^*}{H}\right), \qquad s=\max\{j : \lambda_j \leq 1\}.
$$
Если $z\in Z_0$, то $\lambda_1 \leq 1$ и число $s$ определено корректно. Заметим также, что $\nu(z)\geq1$.

В случае $s\leq r$ получаем
$$\prodd_{i=1}^{2r} \max\left(1, \frac{1}{\lambda_i}\right) = \prodd_{i=1}^{s} (\lambda_i)^{-1} \leq (\lambda_1)^{-r} \leq \nu^{-r}H^{r}.
$$
Если же $s>r$, то ввиду (\ref{7}) и (\ref{9}) имеем
\begin{multline*}
\prodd_{i=1}^{2r}\max\left(1,\frac{1}{\lambda_i}\right) = \prodd_{i=1}^{s} (\lambda_i)^{-1} \leq (\lambda_{2r})^{2r-s}\prodd_{i=1}^{2r} (\lambda_i)^{-1}
\leq r^{O(r)} (\lambda_1^*)^{s-2r} \prodd_{i=1}^{2r} (\lambda_i)^{-1} \leq  \\
r^{O(r)}\left(\frac{p}{H\nu}\right)^{2r-s}\frac{(2r)!H^{2r}}{p^r} = r^{O(r)} \nu^{-r}H^r\left(\frac{H\nu}{p}\right)^{s-r}.
\end{multline*}
Но
$$ \frac{H\nu}{p} \leq \frac{H^2\lambda_1}{p}\leq\frac{H^2}{p}\leq1,
$$
и, значит, в обоих случаях справедливо
$$ \prodd_{i=1}^{2r} \max\left(1,\frac{1}{\lambda_i}\right) \leq r^{O(r)}\nu^{-r}H^r.
$$
Учитывая (\ref{8}), получаем
\begin{equation}
f_0(z) \leq r^{O(r)}\nu(z)^{-r}H^r. \label{10}
\end{equation}
Положим $J:=[\log_2H]+1$ и разобьем множество $Z_0$ на подмножества $Z_j, \quad j=1,\ldots,J$, где
$$ Z_j=\{z\in Z_0 : 2^{j-1} \leq \nu(z) < 2^j \}.
$$
Заметим, что $Z_j \subset Z^1_j \cup Z^2_j$, где
$$ Z^1_j=\{z\in Z_0 : \lambda_1(z)H < 2^j \}, \qquad Z^2_j=\{z\in Z_0 : p\lambda_1^*(z)/H < 2^j\}.
$$
Если $z\in Z^1_j$, то найдется ненулевой вектор $v=(v_1,\ldots,v_{2r})\in\Gamma_z$ такой, что $\max_i|v_i|<2^j$. Вектор $v$ однозначно определяет $z$, так как
$$z = \summ_{i=r+1}^{2r} v_ia_{i-r}\left(\summ_{i=1}^r v_ia_i \right)^{-1}.
$$
(Случай, когда оба вектора
$$ \summ_{i=r+1}^{2r} v_ia_{i-r} \quad \mbox{и} \quad \summ_{i=1}^r v_ia_i
$$
равны нулю в $\F_{p^r}$ невозможен, так как $\max_i|v_i|<2^J\leq p$). Поэтому $\left| Z_1^j\right| <(2^{j+1})^{2r}$. Аналогично, для $z\in Z^2_j$ найдется ненулевой вектор $u=(u_1,\ldots,u_{2r})\in p\Gamma^{*}_z$ такой, что $\sum_i |u_i| <2^j$. Вектор $u$ однозначно определяет $z$. В самом деле, предположим противное. Тогда найдутся различные $z'$ и $z''$ такие, что $u\in p\Gamma^*_{z'}$ и $u\in p\Gamma^*_{z''}$. Возьмем произвольный элемент $x\in\F_{p^r}$, $x=\summ_{i=1}^r x_ia_i$. Пусть
$$ y'=xz'=\summ_{i=1}^{r} y'_ia_i, \quad\quad  y''=xz''=\summ_{i=1}^{r} y''_ia_i .
$$
Тогда $(x_1,\ldots,x_r,y_1',\ldots,y_r')\in\Gamma_{z'}$, \quad $(x_1,\ldots,x_r,y_1'',\ldots,y_r'')\in\Gamma_{z''}$, и мы имеем
\begin{equation}
\summ_{i=1}^{r} y_iu_{i+r}  \equiv 0 \pmod{p},   \label{11}
\end{equation}
где $y_i=y'_i-y''_i$. Так как уравнение $xz'-xz''=y$ разрешимо (относительно $x$) для любого $y\in\F_{p^r}$, то сравнение (\ref{11}) выполнено при всех $(y_1,\ldots,y_r)\in\mathbb{Z}^r$. Значит,
$$ u_{r+1} \equiv \ldots \equiv u_{2r} \equiv 0 \pmod{p},
$$
и в силу того, что $\max_i|u_i|<2^J\leq p$, мы получаем
$$ u_{r+1}=\ldots=u_{2r}=0.
$$
Отсюда следует, что
$$ \summ_{i=1}^r x_iu_i \equiv 0 \pmod{p}
$$
при всех $(x_1,\ldots,x_r)\in\mathbb{Z}^r$, и, стало быть, $u_1=\ldots=u_{2r} =0$, что противоречит выбору вектора $u$. Итак, значит, вектор $u$ однозначно определяет $z$, а поэтому $|Z^2_j|<(2^{j+1})^{2r}$.

Объединяя полученные выше оценки, имеем
$$ |Z_j| \leq |Z_j^1|+|Z_j^2| < 2^{2r+1}2^{2jr},
$$
$$ \summ_{z\in Z_j} f_0(z)^2 \leq r^{O(r)}(2^{j-1})^{-2r}H^{2r}|Z_j| \leq r^{O(r)}H^{2r}
$$
$$\summ_{z\in Z_0} f_0(z)^2 =
\summ_{j=1}^{J}\summ_{z\in Z_j} f_0(z)^2 \leq r^{O(r)} H^{2r}\log p.
$$
Применяя (\ref{4}), получаем утверждение леммы.

\bigskip
\bigskip

Для натурального $H\leq p/2$ и мультипликативного характера $\chi$ на $\F_{p^r}$ положим
$$ \Delta(H,\chi) = \max\limits_B \left|\summ_{x\in B} \chi(x) \right| |B|^{-1},
$$
где максимум берется по всем параллелепипедам вида (\ref{B}) таким, что
\begin{equation}
H\leq H_i \leq 2H, \quad\quad i=1,\ldots,r. \label{12}
\end{equation}
Иногда мы будем писать $\Delta(H)$ вместо $\Delta(H.\,\chi)$ для фиксированного характера $\chi$. Заметим, что если рёбра параллепипеда $B$ удовлетворяют более слабым, чем (\ref{12}), условиям, а именно,
$$ H \leq H_i \leq p, \quad\quad i=1,\ldots,r,
$$
то его можно разбить на параллелепипеды, рёбра которых удовлетворяют условиям (\ref{12}). Поэтому
$$ \left|\summ_{x\in B} \chi(x) \right| \leq \Delta(H,\chi)|B|,
$$
и для того, чтобы доказать теорему C, достаточно проверить, что при всех $0<\varepsilon\leq1/4$ справедливо
\begin{equation}
\Delta([p^{1/4+\varepsilon}],\chi) \ll \frac{1}{\varepsilon}r^{O(1)}p^{-\varepsilon^2/2}. \label{13}
\end{equation}

\bigskip

Приведем еще один вспомогательный результат.

\textsc{Лемма 3} [2, лемма 2]. \textit{Пусть множество $B$ определено соотношением (\ref{B}), где $N_i,H_i>0$ --- целые числа, удовлеторяющие (\ref{12}). Предположим, что $H,\tilde{H}$ --- натуральные числа, $\tilde{H}\leq H/2$, и элемент $u=\sum_{i=1}^r u_ia_i$ удовлетворяет условию}
$$ 1 \leq u_i \leq \tilde{H}, \quad\quad i=1,\ldots,r.
$$
\textit{Пусть, далее, $\chi$ --- мультипликативный характер в $\F_{p^r}$. Тогда}
$$ \left| \summ_{x\in B} \chi(x) - \summ_{x\in B} \chi(x+u) \right| \leq 6r\Delta(\tilde{H},\chi)|B|\frac{\tilde{H}}{H}.
$$
\smallskip

Следующая лемма основана на хорошо известном методе, разработанном Бёрджесом (см., например, [9],[10],[11]), который позволяет получить оценку сумм характеров через лемму 2 и лемму E.

\smallskip

\textsc{Лемма 4.} \textit{Пусть $s,H,\tilde{H}$ --- натуральные числа, и}
$$ G:=[p^{r/(2s)}] \leq \tilde{H} \leq \frac{H}{2}, \quad\quad H\leq p^{1/2}.
$$
\textit{Пусть также $\chi$ --- нетривиальный мультипликативный характер на $\F_{p^r}$. Тогда имеет место следующее неравенство:}
$$ \Delta(H,\chi) \leq 6r\Delta(\tilde{H})\frac{\tilde{H}}{H}+ O(sr^{O(r/s)}(\log p)^{1/(2s)}(H\tilde{H}G^{-1}p^{-1/2})^{-r/(2s)}).
$$
\textit{Док-во}. Положим $I=[1,G]\cap \mathbb{Z}$, $H_0=[\tilde{H}/G]$, и
$$B'_0 = \left\{\summ_{i=1}^r x_ia_i : 1\leq x_i \leq H_0, \, i=1,\ldots,r \right\} .
$$
Заметим, что любой элемент $x=yz,$ $y\in B'_0$, $z\in I$, может быть представлен в виде
$$ x = \summ_{i=1}^r x_ia_i : \quad\quad 1\leq x_i \leq \tilde{H}, i=1,\ldots,r.
$$
Согласно лемме 3, для всех $y\in B'_0$, $z\in I$ справедливо
$$ \left| \summ_{x\in B} \chi(x+yz) - \summ_{x\in B} \chi(x) \right| \leq 6r\Delta(\tilde{H},\chi)|B|\frac{\tilde{H}}{H} .
$$
Поэтому
\begin{equation}
\left| \summ_{x\in B} \chi(x) - \frac{1}{|B_0'|G}\summ_{x\in B,y\in B'_0,z\in I} \chi(x+yz) \right|  \leq 6r\Delta(\tilde{H},\chi)|B|\frac{\tilde{H}}{H} .
\label{17}
\end{equation}
Далее,
\begin{multline*}
\left| \summ_{x\in B,y\in B'_0,z\in I} \chi(x+yz) \right| \leq \summ_{x\in B,y\in B'_0} \left|\summ_{z\in I}\chi(x+yz)\right| \\
= \summ_{x\in B,y\in B'_0} \left|\summ_{z\in I}\chi(xy^{-1}+z)\right| = \summ_{u\in\F_{p^r}} \omega(u)\left|\summ_{z\in I} \chi(u+z)\right|,
\end{multline*}
где
$$ \omega(u)=\left| \left\{(x,y)\in B \times B'_0 : \frac xy =u\right\} \right| .
$$
Дважды используя неравенство Гёльдера, получаем
\begin{multline}
\left|\summ_{x\in B,y\in B'_0,z\in I} \chi(x+yz)\right| \leq \\
\left(\summ_{u\in\F_{p^r}} \omega(u)\right)^{1-1/s}\left(\summ_{u\in\F_{p^r}} \omega(u)^2\right)^{1/(2s)}
\left(\summ_{u\in\F_{p^r}} \left|\summ_{z\in I} \chi(u+z)\right|^{2s} \right)^{1/(2s)}.
\label{18}
\end{multline}
Очевидно,
\begin{equation}
\summ_{u\in\F_{p^r}} \omega(u) = |B||B'_0|. \label{19}
\end{equation}
Для того, чтобы оценить сумму $\summ_{u\in\F_{p^r}} \omega(u)^2$, введём следующее обозначение: для множества $A\subset \F_{p^r}$ через $E(A)$ будем обозначать число решений уравнения
$$ x^1x^2=x^3x^4, \quad\quad x^1,x^2,x^3,x^4 \in A.
$$
Положим $B^*=B\setminus\{0\}$. Имеем
\begin{gather*} \summ_{u\in\F_{p^r}} \omega(u)^2 = \omega(0)^2+\summ_{u\in\F^{*}_{p^r}} \omega(u)^2 \\
= \omega(0)^2+\left| \{ (x_1,x_2,y_1,y_2)\in B^*\times B^* \times B'_0 \times B'_0 : x_1y_2 = x_2y_1 \}\right| \\
\leq |B_0'|^2+ \summ_{\nu\in\F^{*}_{p^r}} \left| \left\{ (x_1,x_2)\in B^*\times B^* : \frac{x_1}{x_2}=\nu \right\} \right| \left| \left\{ B'_0 \times B'_0 : \frac{y_1}{y_2} = \nu \right\}\right| \\
\leq |B_0'|^2+E(B)^{1/2}E(B_0')^{1/2} .
\end{gather*}
Используя лемму 2, и принимая во внимание, что $|B_0'|\leq |B|$, получим
\begin{equation}
\summ_{u\in\F_{p^r}} \omega(u)^2 \leq r^{O(r)}|B||B'_0|\log p .
\end{equation}
Далее, мы используем оценку последней суммы в (\ref{18}), полученную в [11]. Для полноты мы воспроизводим здесь доказательство.
\begin{multline*}\summ_{u\in\F_{p^r}} \left|\summ_{z\in I} \chi(u+z)\right|^{2s} = \\
\summ_{z_1,\ldots,z_{2s}\in I}\left|\summ_{u\in\F_q} \chi\left((u+z_1)\ldots u(u+z_s)(u+z_{s+1})^{q-2}\ldots(u+z_{2s})^{q-2}\right)\right|.
\end{multline*}
Последняя сумма уже оценивалась нами в частном случае, когда $\chi$ -- квадратичный характер, при доказательстве леммы 1. В общем случае работают аналогичные рассуждения.
Мы можем применить лемму E для тех наборов $(z_1,\ldots,z_{2s})$, в которых хотя бы один элемент встретился ровно один раз. Для таких наборов получим (см. оценку на сумму )
$$ \left|\summ_{u\in\F_q} \chi\left((u+z_1)\ldots u(u+z_s)(u+z_{s+1})^{q-2}\ldots(u+z_{2s})^{q-2}\right)\right| < 4sp^{r/2}.
$$
Сумму по наборам $(z_1,\ldots,z_{2s})$, в которых каждый элемент встретился как минимум дважды, оценим произведением $p^r$ на количество таких наборов. Так как все элементы $z_1,\ldots,z_{2s}$ лежат в некотором подмножестве множества $I$, содержащем $\min(s,G)$ элементов, и для конкретного такого подмножества есть не более $s$ способов выбрать каждый элемент, то количество интересующих нас наборов не превосходит $G^ss^{2s}$. Поэтому
$$ \summ_{u\in\F_{p^r}} \left|\summ_{z\in I} \chi(u+z)\right|^{2s} \leq G^ss^{2s}p^r + 4sG^{2s}p^{r/2},
$$
$$ \left(\summ_{u\in\F_{p^r}} \left|\summ_{z\in I} \chi(u+z)\right|^{2s}\right)^{1/(2s)} \leq sG^{1/2}p^{r/(2s)} + 2Gp^{r/(4s)}.
$$
По определению числа $G$ имеем
\begin{equation}
\left(\summ_{u\in\F_{p^r}} \left|\summ_{z\in I} \chi(u+z)\right|^{2s}\right)^{1/(2s)} \leq (2s+2)Gp^{r/(4s)}. \label{21}
\end{equation}
Подставляя (\ref{19})--(\ref{21}) в (\ref{18}), получаем
\begin{multline}
\left|\summ_{x\in B,y\in B'_0,z\in I} \chi(x+yz)\right| \ll  sr^{O(r/s)}\left(|B||B'_0|\right)^{1-1/(2s)} (\log p)^{1/(2s)} Gp^{r/(4s)} \\
= sr^{O(r/s)}(\log p)^{1/(2s)}|B||B'_0|G \left(|B||B'_0|\right)^{-1/(2s)}p^{r/(4s)} \\
\leq sr^{O(r/s)} (\log p)^{1/(2s)}|B||B'_0|G\left( H\left[\frac{\tilde{H}}{G}\right]p^{-1/2}\right)^{-r/(2s)} \\
\leq 2^{r/(2s)}sr^{O(r/s)}(\log p)^{1/(2s)}|B||B'_0|G\left(H\tilde{H}G^{-1}p^{-1/2} \right)^{-r/(2s)}.
\end{multline}
Из последнего неравенства и (\ref{17}) следует, что
$$ \left|\summ_{x \in B} \chi(x)\right| \leq 6r\Delta(\tilde{H},\chi)|B|\frac{\tilde{H}}{H} + O(sr^{O(r/s)}(\log p)^{1/(2s)}|B|\left(H\tilde{H}G^{-1}p^{-1/2} \right)^{-r/(2s)}).
$$
Так как выбор $B$ был произвольным, то лемма доказана.

\bigskip
\bigskip

Предположим в лемме 4, что
\begin{equation} s \geq 2r.  \label{22}
\end{equation}
Введем параметр $\alpha$, удовлетворяющий условию
\begin{equation} 0< \alpha \leq \frac{1}{(12r)^2}. \label{23}
\end{equation}
Предположим также, что
\begin{equation} 0 < \varepsilon \leq \frac14, \quad\quad H=[p^{1/4+\varepsilon}]. \label{24}
\end{equation}
Положим
$$ H_i = [\alpha^iH], \quad\quad i\geq0.
$$
Пусть $I$ -- наибольшее целое число такое, что $H_i\geq G$, где $G$ определено в лемме 4. Так как $H\geq G$ в силу (\ref{22}) и (\ref{24}), то число $I$ определено корректно.

Применяя лемму 4 к $H_i$ и $H_{i+1}$, мы можем последовательно оценить $\Delta(H_i)$ через $\Delta(H_{i+1})$, $i=0,\ldots,I-1$ :
$$ \Delta(H_i) \leq 6r\Delta(H_{i+1})\frac{H_{i+1}}{H_i} + O\left(sr^{O(r/s)}(\log p)^{1/(2s)}(H_iH_{i+1}G^{-1}p^{-1/2})^{-r/(2s)}\right).
$$
Поэтому
\begin{multline}
\Delta(H) \leq (6r)^I\Delta(H_I)\frac{H_I}{H} + \\
sr^{O(r/s)}(\log p)^{1/(2s)}(Gp^{1/2})^{r/(2s)}\summ_{i=0}^{I-1} O\left((6r)^i\frac{H_i}{H}(H_iH_{i+1}) ^{-r/(2s)}\right). \label{25}
\end{multline}
Далее,
\begin{multline}
(6r)^I\Delta(H_I)\frac{H_I}{H} \leq (6r)^{\log H / (-\log \alpha)}\frac{H_I}{H} \\
\leq (6r)^{\log p / (-2\log \alpha)}\alpha^{-1}\frac{2G}{H} = 2p^{\log (6r) / (-2\log \alpha)}\alpha^{-1}\frac{G}{H} \quad . \label{26}
\end{multline}
Принимая во внимание (\ref{22}) и (\ref{23}), получим
\begin{multline}
\summ_{i=0}^{I-1} (6r)^i\frac{H_i}{H}(H_iH_{i+1})^{-r/(2s)} \leq \summ_{i=0}^{I-1} (6r)^i\alpha^{i}\left(\frac{H^2\alpha^{2i+1}}{4}\right)^{-r/(2s)} \\
 \leq 2H^{-r/s}\alpha^{-r/(2s)}\summ_{i=0}^{I-1} (6r)^i\alpha^{i(1-r/s)} \leq 2H^{-r/s}\alpha^{-r/(2s)}\summ_{i=0}^{I-1} (6r)^i\alpha^{i/2} \\
 \leq 2H^{-r/s}\alpha^{-r/(2s)}\summ_{i=0}^{I-1}2^{-i} \leq 4H^{-r/s}\alpha^{-r/(2s)} . \label{27}
\end{multline}
Подставляя (\ref{26}) и (\ref{27}) в (\ref{25}), находим
\begin{multline}
\Delta(H) \leq 2p^{\log 6r / (-2\log \alpha)}\alpha^{-1}\frac{G}{H} + O\left(sr^{O(r/s)}(\log p)^{1/(2s)}(Gp^{1/2})^{r/(2s)} H^{-r/s}\alpha^{-r/(2s)}\right) \\
\leq 4p^{\log 6r / (-2\log \alpha)}\alpha^{-1}p^{r/(2s)-1/4-\varepsilon}+O\left(sr^{O(r/s)}(\log p)^{1/(2s)}\alpha^{-1}p^{-(r/s)(\varepsilon - r/(4s))}\right). \label{28}
\end{multline}
Для завершения доказательства теоремы С выберем параметры $s$ и $\alpha$ следующим образом:
$$ s=\left[\frac{r}{2\varepsilon}\right]+1, \quad\quad \alpha=(6r)^{-3} .
$$
Легко видеть, что при этом условия (\ref{22}) и (\ref{23}) выполнены. Так как
$$ \frac{4\varepsilon}{3} \leq \frac {r}{s} \leq 2\varepsilon,
$$
то
$$ \frac{r}{s}\left(\varepsilon - \frac{r}{4s}\right) \geq \frac23\varepsilon^2 \, ;
$$
следовательно,
\begin{multline*}(\log p)^{1/(2s)}\alpha^{-1}p^{-(r/s)(\varepsilon - r/(4s))} \leq (6r)^3(\log p)^{\varepsilon/r}p^{-2\varepsilon^2/3} \\
= (6r)^3\exp\left(\frac{\varepsilon}{r}\log\log p - \frac{\varepsilon^2}{6}\log p\right)p^{-\frac{\varepsilon^2}{2}} \ll r^3p^{-\frac{\varepsilon^2}{2}} .
\end{multline*}
Наконец,
$$ p^{\log 6r / (-2\log \alpha)}\alpha^{-1}p^{r/(2s)-1/4-\varepsilon} \leq (6r)^3 p^{1/6} p^{-1/4} \ll r^3p^{-1/12} .
$$
Значит, неравенство (\ref{28}) дает нам
$$ \Delta(H) \ll r^3p^{-1/12} + \frac{r}{\varepsilon}r^{3+O(\varepsilon)}p^{-\varepsilon^2/2} \leq \frac{r^{4+O(\varepsilon)}p^{-\varepsilon^2/2}}{\varepsilon}.
$$
Таким образом, неравенство  (\ref{28}) влечет неравенство  (\ref{13}), и теорема С доказана.

\end{fulltext}

\end{document}